\documentstyle[12pt]{article}

\makeatletter
\newdimen\normalarrayskip
\newdimen\minarrayskip
\normalarrayskip\baselineskip
\minarrayskip\jot
\newif\ifold \oldtrue 
\def\arraymode{\ifold\relax\else\displaystyle\fi}
\def\eqnumphantom{\phantom{(\theequation)}}
\def\@arrayskip{\ifold\baselineskip\z@\lineskip\z@
\else
\baselineskip\minarrayskip\lineskip2\minarrayskip\fi}
\def\@arrayclassz{\ifcase \@lastchclass \@acolampacol \or
\@ampacol \or \or \or \@addamp \or
\@acolampacol \or \@firstampfalse \@acol \fi
\edef\@preamble{\@preamble
\ifcase \@chnum
\hfil$\relax\arraymode\@sharp$\hfil
\or $\relax\arraymode\@sharp$\hfil
\or \hfil$\relax\arraymode\@sharp$\fi}}
\def\@array[#1]#2{\setbox\@arstrutbox=\hbox{\vrule
height\arraystretch \ht\strutbox
depth\arraystretch \dp\strutbox
width\z@}\@mkpream{#2}\edef\@preamble{\halign
\noexpand\@halignto
\bgroup \tabskip\z@ \@arstrut \@preamble \tabskip\z@ \cr}%
\let\@startpbox\@@startpbox \let\@endpbox\@@endpbox
\if #1t\vtop \else \if#1b\vbox \else \vcenter \fi\fi
\bgroup \let\par\relax
\let\@sharp##\let\protect\relax
\@arrayskip\@preamble}
\def\eqnarray{\stepcounter{equation}%
\let\@currentlabel=\theequation
\global\@eqnswtrue
\global\@eqcnt\z@
\tabskip\@centering
\let\\=\@eqncr
$$%
\halign to \displaywidth\bgroup
\eqnumphantom\@eqnsel\hskip\@centering
$\displaystyle \tabskip\z@ {##}$%
\global\@eqcnt\@ne \hskip 2\arraycolsep
$\displaystyle\arraymode{##}$\hfil
\global\@eqcnt\tw@ \hskip 2\arraycolsep
$\displaystyle\tabskip\z@{##}$\hfil
\tabskip\@centering
&{##}\tabskip\z@\cr}
\begingroup\ifx\undefined\newsymbol \else\def\input#1 {\endgroup}\fi

\def\nl#1#2{\mathop{#2}\limits_{#1}}

\newcommand{\pl}{\partial}

\newcommand{\bqq}{\begin{equation} \label}
\newcommand{\eeq}{\end{equation}}

\begin{document}

\large
\centerline{\Large{\bf{Rigid 6-dimensional $h$-spaces of constant 
curvature
\footnote{The research was partially supported by the RFBR
grant 01-02-17682-a and by the INTAS grant 00-00334.}}}}

\bigskip

\centerline{{\large
Zolfira Zakirova}\footnote{Kazan State University,
e-mail: Zolfira.Zakirova@soros.ksu.ru, Zolfira.Zakirova@ksu.ru}}

\bigskip

\abstract{\small  In this paper, we continue studying the $6$-dimensional
pseudo-Riemannian space $V^6(g_{ij})$
with signature $[++----]$, which admits
projective motions, i. e. continuous transformation groups preserving
geodesics. In particular, we determine a necessary and sufficient
condition that the 6-dimensional rigid $h$-spaces have constant curvature.
}

\begin{center}
\rule{5cm}{1pt}
\end{center}

\bigskip
\setcounter{footnote}{0}
\section{Introduction}

The general method of determining pseudo-Riemannian manifolds that
admit some nonhomothetic projective group $G_r$ has been developed
by A.V.Aminova \cite{am2}. A.V.Aminova has classified
all the Lorentzian manifolds
of dimension $\geq 3$ that admit nonhomothetic projective or affine
infinitesimal transformations. In each case, there were determined the
corresponding maximal projective and affine Lie algebras.
This problem is not solved for pseudo-Riemannian spaces with arbitrary 
signature.
In the series of works, we investigate the $6$-dimensional
pseudo-Riemannian space $V^6(g_{ij})$
with signature $[++----]$.

Remind that in order to find a pseudo-Riemannian space admitting a 
nonhomothetic
infinitesimal projective transformation, one needs to integrate the
Eisenhart equation \cite{ezen1}
\bqq{1}
h_{ij,k}=2g_{ij} \varphi_{,k}+g_{ik} \varphi_{,j}+
g_{jk} \varphi_{,i}.
\eeq
Pseudo-Riemannian manifolds for which there exist nontrivial solutions
$h_{ij}\ne cg_{ij}$ to the Eisenhart equation are called {\it $h$-spaces}.
It is known that the problem of determining such spaces depends on the 
type of
the $h$-space, i. e. on the type of the bilinear form
$L_{X}g_{ij}$ determined by the characteristic of the $\lambda$-matrix
$( h_{ij}-\lambda g_{ij})$. The number of possible types depends on the
dimension and the signature of the $h$-space. In our case of
six-dimensional spaces $V^6$ with signature $[++----]$, there is a variety
of possible types. We restrict it considering the rigid $h$-spaces:
$h$-spaces with distinct bases of prime divisors of the
$\lambda$-matrix are called rigid $h$-spaces \cite{pet,zak}.
In the papers \cite{zak1}-\cite{zak4}
(see also \cite{zak5}), we found 6-dimensional rigid
$h$-spaces of the $[2211]$, $[321]$, $[33]$, $[411]$ and $[51]$ types.

For further study we need a necessary and sufficient
condition of constant curvature of this $h$-space.
This condition is known to be expressed by the formula
\bqq{4.1}
R_{jkl}^i=K({\delta_k}^i g_{jl}-{\delta_l}^i g_{jk}),
\quad
K={\rm const}.
\eeq
Since the projective-group properties of the rigid h-spaces of the $[21111]$
and $[3111]$
types have been investigated by A. V. Aminova, we do not discuss them here.

\section{The $h$-space of the $[2211]$ type}

Let us investigate the $h$-space of the $[2211]$ type.
In this case, metric has the following form
\bqq{2.28}
g_{ij}dx^idx^j=e_2 (f_4-f_2)^2  \Pi_{\sigma} (f_{\sigma}-f_2)
\lbrace 2A  dx^1dx^2-A^2 \Sigma_1 (dx^2)^2 \rbrace+
\eeq
$$
e_4 (f_2-f_4)^2  \Pi_{\sigma} (f_{\sigma}-f_4) \lbrace 2 \tilde{A}
 dx^3dx^4-\tilde{A}^2 \Sigma_2 (dx^4)^2 \rbrace+
\sum_{\sigma} e_{\sigma} {\Pi'_i}(f_i-f_{\sigma}) (dx^{\sigma})^2,
$$
where
\bqq{2.28.4}
A=\epsilon x^1+\theta(x^2),
\quad
\tilde A=\tilde{\epsilon} x^3+\omega(x^4),
\eeq
$$
\Sigma_1=2 (f_4-f_2)^{-1}+ \sum_{\sigma}(f_{\sigma}-f_2)^{-1},
\quad
\Sigma_2=2 (f_2-f_4)^{-1}+ \sum_{\sigma}(f_{\sigma}-f_4)^{-1},
$$
$f_2=\epsilon x^2$, $f_4=\tilde{\epsilon}x^4+a$,
$\epsilon, \tilde{\epsilon}=0, 1$,
$a$ is a constant nonzero when $\tilde{\epsilon}=0$,
$f_{\sigma}(x^{\sigma})$, $\theta(x^2)$,  $\omega(x^4)$ are
arbitrary functions, $e_i= \pm 1$.

Calculating components of the curvature tensor of our $h$-space,
we observe that the nonzero components are
\bqq{4.2}
R_{b_p a_p c_p}^{a_p}=\chi_p g_{b_p c_p},
\eeq
$$
R_{b_p \sigma c_p}^{\sigma}=\rho_{\sigma p} g_{a_p b_p}-
\frac{\chi_p-\rho_{\sigma p}}{f_{\sigma}-f_p} A_p g_{a_p c_p} 
\delta_{b_p}^p,
\; (c_p \ne b_p),
$$
$$
R_{\sigma b_p \sigma}^{a_p}=g_{\sigma \sigma}\lbrace \rho_{\sigma p}
\delta_{b_p}^{a_p} -\frac{\chi_p-\rho_{\sigma p}}{f_{\sigma}-f_p}A_p
\delta_{b_p}^p \delta_{c_p}^{a_p}\rbrace, \; (c_p \ne b_p),
$$
$$
R_{\sigma \tau \sigma}^{\tau}=g_{\sigma \sigma}\frac{\rho_{\tau p}
(f_{\tau}-f_p)-\rho_{\sigma p}(f_{\sigma}-f_p)}{f_{\tau}-f_{\sigma}},
\; (\sigma \ne \tau),
$$
$$
R_{a_q b_p b_q}^{a_p}=\lbrace \rho_{pq} g_{a_q b_q} -\frac{\chi_q-\rho_{pq}}
{f_p-f_q}A_q g_{a_q c_q} \delta_{b_q}^q\rbrace \delta_{b_p}^{a_p}-
$$
$$
\lbrace \frac{\chi_p-\rho_{pq}}
{f_q-f_p} A_p g_{a_q b_q}+\frac{ {\nl{l=2,4}\sum} (\chi_{l}-\rho_{pq}) }
{(f_q-f_p)^2}A_p A_q g_{a_q c_q} \delta_{b_q}^q \rbrace \delta_{b_p}^{p}
\delta_{c_p}^{a_p},
$$
$$
(c_p \ne b_p, c_q \ne b_q, p \ne q),
$$
where
\bqq{4.3}
\chi_p=B_p+\rho_p,
\quad
B_2=\frac{\epsilon \theta'}{A^2 g_{12}},
\;
B_4=\frac{\tilde{\epsilon} \omega'}{\tilde {A}^2 g_{34}},
\quad
A_2=A,
\;
A_4=\tilde A,
\eeq
\bqq{4.4}
\rho_p=-\frac{1}{4} \sum_{\sigma} \frac{(f'_{\sigma})^2}{(f_{\sigma}-f_p)^2
g_{\sigma \sigma}},
\quad
\rho_{pq}=-\frac{1}{4} \sum_{\sigma} \frac{(f'_{\sigma})^2}
{(f_{\sigma}-f_p)(f_{\sigma}-f_q)g_{\sigma \sigma}},
\eeq
$$
\rho_{\sigma p}=-\frac{1}{4}\frac{(f'_{\sigma})^2}{(f_{\sigma}-f_p)
g_{\sigma \sigma}} \lbrace \frac{2 {f''}_{\sigma}}{(f'_{\sigma})^2}-
\frac{1}{f_{\sigma}-f_p}+
\sum_{i, i \ne \sigma} (f_i-f_{\sigma})^{-1} \rbrace-
$$
$$
\frac{1}{4} \sum_{\gamma, \gamma \ne \sigma} \frac{(f'_{\gamma})^2}
{(f_{\gamma}-f_p)(f_{\gamma}-f_{\sigma})g_{\gamma \gamma}}.
$$
Here $p, q=2,4$, $\sigma=5, 6$, the indices $a_2$  are equal to either 1 
or 2,
while the indices $a_4$ are 3 or 4.

In particular, from (\ref{4.1}) it follows that
$$
R_{112}^1=R_{1\sigma 2}^{\sigma},
\quad
R_{334}^3=R_{3\sigma 4}^{\sigma},
\quad
R_{112}^1=R_{142}^4,
\quad
R_{334}^3=R_{324}^2.
$$
Substituting components of the curvature tensor into this equality, one
obtains
$$
(\chi_2-{\rho}_{\sigma 2})g_{12}=0,
\quad
(\chi_4-{\rho}_{\sigma 4})g_{34}=0,
\quad
(\chi_2-{\rho}_{24})g_{12}=0,
\quad
(\chi_4-{\rho}_{24})g_{34}=0,
$$
since $g_{12} \ne 0$ and $g_{34} \ne 0$, then
\bqq{4.5}
\chi_2-{\rho}_{\sigma 2}=0,
\quad
\chi_4-{\rho}_{\sigma 4}=0,
\eeq
$$
\chi_2-{\rho}_{24}=0,
\quad
\chi_4-{\rho}_{24}=0.
$$
Let us prove that the condition  $\chi_2-{\rho}_{\sigma 2}=0$ is equivalent
to the conditions
$$
\rho_{2}-\rho_{\sigma 2}=\epsilon=0.
$$
Indeed, differentiating the equality
$$
\chi_2-{\rho}_{\sigma 2}=\rho_{2}-\rho_{\sigma 2}+\frac{\epsilon \theta'}
{A^2 g_{12}}=0
$$
with respect to $x^1$, we obtain $\epsilon=0$ and, therefore,
$\rho_{2}-\rho_{\sigma 2}=0$. Inversely, if $\rho_{2}-\rho_{\sigma 2}=
\epsilon=0$, then $\chi_2-{\rho}_{\sigma 2}=0$.
Similarly we can find equivalent conditions for
other relations (\ref{4.5}).

Let us prove that the necessary conditions of constant curvature of the
$h$-space of the $[2211]$ type,
\bqq{4.6}
\chi_p-{\rho}_{\sigma p}=\chi_p-{\rho}_{24}=0,
\eeq
or equivalent conditions
\bqq{4.7}
\rho_p-{\rho}_{\sigma p}=\rho_p-{\rho}_{24}=\epsilon=\tilde \epsilon=0,
\eeq
are sufficient. The relations that follows from (\ref{4.4}),
$$
\pl_{\sigma} \rho_p=-f'_{\sigma}\frac{\rho_p-\rho_{\sigma 
p}}{f_{\sigma}-f_p},
\quad
\pl_1 \rho_p=\pl_3 \rho_p=0,
$$
$$
\pl_p \rho_p=-\epsilon_p \sum_{\sigma} 
\frac{{f'_{\sigma}}^2}{(f_{\sigma}-f_p)^3 g_{\sigma \sigma}},
\quad
\pl_p \rho_q=-\frac{ \epsilon_p}{2} \sum_{\sigma} 
\frac{{f'_{\sigma}}^2}{(f_{\sigma}-f_q)^2(f_{\sigma}-f_p) g_{\sigma 
\sigma}},
$$
where $p \ne q$, $\epsilon_2=\epsilon$, $\epsilon_4=\tilde \epsilon$,
along with (\ref{4.7}), lead to the identities
$$
\chi_p=\rho_{\sigma p}=\rho_{24}=\rho_p={\rm const}.
$$
Using these results along with conditions $\epsilon=\tilde\epsilon=0$,
one can prove that components of the curvature tensor (\ref{4.2})
satisfy conditions (\ref{4.1}) with $K=\rho_p={\rm const}$.

\section{The $h$-space of the $[321]$ type}

Metric of this h-space can be written in the form
\bqq{2.29-4}
g_{ij}dx^idx^j=e_3 (f_5-\epsilon x^3)^2(f_6-\epsilon x^3)
\lbrace (dx^2)^2 +4A  dx^1dx^3+
\eeq
$$
2(\epsilon x^1-2 A \Sigma_1)dx^2 dx^3+
((\epsilon x^1)^2-4\epsilon x^1 A \Sigma_1+
4A^2 \Sigma_3) (dx^3)^2 \rbrace+
$$
$$
e_5 (\epsilon x^3-f_5)^3
(f_6-f_5)
\lbrace 2\tilde{A}  dx^4dx^5-
\Sigma_4 {\tilde{A}}^2(dx^5)^2 \rbrace+
e_6(f_5-f_6)^2(f_5-f_6)^3
 (dx^6)^2,
$$
where
\bqq{2.29-5}
A=\epsilon x^2+\theta(x^3),
\quad
\tilde{A}=\tilde{\epsilon} x^4+\omega(x^5),
\eeq
$$
\Sigma_1=(f_6-\epsilon x^3)^{-1}+
2(f_5-\epsilon x^3)^{-1},
\quad
\Sigma_2=(f_6-\epsilon x^3)^{-2}+
2(f_5-\epsilon x^3)^{-2},
$$
$$
2\Sigma_3=(\Sigma_1)^2-\Sigma_2,
\quad
\Sigma_4=3(\epsilon x^3-f_5)^{-1}+
(f_6-f_5)^{-1},
$$
$f_1=f_2=f_3=\epsilon x^3$,
$f_4=f_5=\tilde{\epsilon}x^5+a$, $f_6$ is an arbitrary function of the 
variable
$x^6$, $\epsilon, \tilde{\epsilon}=0, 1$, $\epsilon\ne0$ when
$\tilde\epsilon=0$, and, inversely,
$\tilde\epsilon\ne0$ when $\epsilon=0$,
$a$ is a constant nonequal to zero when $\tilde\epsilon=0$,
$e_3, e_5, e_6 =\pm 1$.

 From (\ref{4.1}) one finds
$$
R_{123}^2= R_{1r3}^r,
\quad
(r~=~5,6).
$$
Substituting into this equality the corresponding components of the
curvature tensor of our $h$-space
$$
R_{123}^2=\chi_3 g_{13},
\quad
R_{153}^5=\rho_{53} g_{13},
\quad
R_{163}^6=\rho_{63} g_{13},
$$
where
\bqq{4.20}
\chi_3=\frac{3{\epsilon}^2}{16 A^2 g_{22}}+\rho_3,
\eeq
$\rho_3$, $\rho_{53}$ and $\rho_{63}$ are defined by formulas 
(\ref{4.4}) with
$p, q=3, 5$ and $\sigma=5, 6$, one gets
$$
(\chi_3-\rho_{r3})g_{13}=0.
$$
Since $g_{13} \ne 0$, then
$$
\chi_3-\rho_{r3}=\frac{3 {\epsilon }^2}
{16 A^2 g_{22}}+{\rho}_{3}-{\rho}_{r3}=0.
$$
Differentiating the last equality with respect to
$x^2$, one obtains $\epsilon=0$ and, therefore, $\rho_3-\rho_{r3}=0$. 
Inversely,
if $\rho_3-\rho_{r3}=\epsilon=0$, then $\chi_3-\rho_{r3}=0$.
Since $f_3 \ne f_5$, from $\rho_3-\rho_{53}=0$ it follows that $f'_6=0$.
Substituting $\epsilon=0$ into the relation
$$
R_{123}^1=\gamma
g_{13}+\frac{3{\epsilon}^2}{8A}\Sigma_1+\frac{3\epsilon}{4A^2}(\theta'-
{\epsilon}^2x^1)=0,
$$
where $R_{123}^1=0$ follows from (\ref{4.1}), one finds
$$
\gamma=-\frac{1}{4} {\nl{\sigma}{\sum}} \frac{(f'_{\sigma})^2}
{(f_{\sigma}-f_p)^3 g_{\sigma \sigma}}=0.
$$
 From (\ref{4.1}), one has
$$
R_{445}^4=R_{4 s 5}^{s},
\quad
(s=3, 6).
$$
Therefore,
$$
(\chi_5-\rho_{s5})g_{45}=0,
$$
since for our $h$-space
$$
R_{445}^4=\chi_5 g_{45}, \quad
R_{435}^3=\rho_{53} g_{45},
\quad
R_{465}^6=\rho_{65} g_{45},
$$
where
\bqq{4.21}
\chi_5=\frac{\tilde\epsilon \theta'}{{\tilde A}^2 g_{45}}+\rho_5.
\eeq
Since $g_{45} \ne 0$, we obtain $\chi_5-\rho_{s5}=0$.
Similarly one can prove that this condition is equivalent to
the conditions $f'_6=\tilde \epsilon=0$.

Let us prove that the necessary conditions of constant curvature of the
$h$-space of the $[321]$ type,
\bqq{4.22}
\chi_p-\rho_{6p}=\chi_p-\rho_{35}=\gamma=0,
\quad
(p=3, 5),
\eeq
or
$$
f'_6=\epsilon=\tilde\epsilon=0,
$$
are sufficient. Note that from (\ref{4.22})
and (\ref{4.4}) $(p,q=3, 5, \sigma=6)$, (\ref{4.20}) and (\ref{4.21})
it follows that
$$
\chi_p=\rho_p=\rho_{6 p}=\rho_{35}=0.
$$
Calculating the other components of the curvature tensor of our 
$h$-space and
substituting them into the last equalities, one gets $h$-space of the
constant curvature.

\section{The $h$-space of the $[33]$ type}

Metric is determined by formulas
\bqq{2.31}
g_{ij}dx^idx^j=e_3 (f_6-f_3)^3 \lbrace (dx^2)^2 +4A  dx^1dx^3+
\eeq
$$
2(\epsilon x^1-2 A \Sigma_1)dx^2 dx^3+((\epsilon x^1)^2-4\epsilon x^1
A \Sigma_1+4A^2 \Sigma_2) (dx^3)^2 \rbrace+
$$
$$
e_6 (f_3-f_6)^3
\lbrace (dx^5)^2+4 \tilde{A}  dx^4dx^6+
2(\tilde{\epsilon} x^4+2 \tilde{A} {\Sigma}_1)dx^5dx^6+
((\tilde{\epsilon} x^4)^2+4\tilde{\epsilon} x^4 \tilde{A} {\Sigma}_1+
4 \tilde{A}^2{\Sigma}_2
( (dx^6)^2 \rbrace,
$$
where
\bqq{2.35}
A=\epsilon x^2+\theta(x^3),
\quad
\tilde{A}=\tilde{\epsilon} x^4+\omega(x^6),
\eeq
$$
\Sigma_1=3 (f_6-f_3)^{-1},
\quad
\Sigma_2=3 (f_6-f_3)^{-2},
$$
$f_3=\epsilon x^3$, $f_6=\tilde{\epsilon}x^6+a$,
$\epsilon, \tilde{\epsilon}=0, 1$, $c$ is a constant, $a$ is a constant 
nonequal
to zero when
$\tilde{\epsilon}=0$, $i_1, j_1=1,2, 3$, $i_2, j_2=4, 5, 6$,  $e_i =\pm 1$,
$\theta(x^3)$, $\omega(x^6)$ are arbitrary functions,
nonzero when $\epsilon=0$ and $\tilde\epsilon=0$ respectively.

Nonzero components
$R_{\nu b_p \mu}^{a_p}$ $(b_p \ne \mu, a_p < b_p)$, $R_{5 a_1
6}^{a_1}$, $R_{6 a_1 5}^{a_1}$, $R_{6 a_1 6}^{a_1}$, $R_{2 a_2 3}^{a_2}$,
$R_{3 a_2 2}^{a_2}$, $R_{3 a_2 3}^{a_2}$ of the curvature tensor
($p=1, 2$, $a_1, b_1=1, 2, 3$, $a_2,
b_2=4, 5, 6$, $\nu$ and $\mu$ are the indices of nonzero
components of metric $g_{\nu \mu}$ (\ref{2.31})) are proportional to 
$\epsilon$
or to $\tilde \epsilon$.
In particular,
$$
R_{123}^2=\frac{3 {\epsilon}^2}{8 A},
R_{456}^5=\frac{3
{\tilde\epsilon}^2}{8 \tilde  A}.
$$
The equalities $R_{123}^2=R_{163}^6$ and
$R_{456}^5=R_{436}^3$ obtained from (\ref{4.1}) lead to
$\epsilon=\tilde\epsilon=0$.  Obviously, when $\epsilon=\tilde\epsilon=0$,
all the
components of the curvature tensor of the $h$-space we consider are equal
to zero.

\section{The $h$-space of the $[411]$ type}

Metric of this space is the following
\bqq{2.45}
g_{ij}dx^idx^j=e_4{\Pi}_{\sigma}(f_\sigma-\epsilon x^4)
\lbrace 6 Adx^1dx^4+2 dx^2dx^3+
\eeq
$$
2(2\epsilon x^2-3A\Sigma_1)dx^2dx^4 -\Sigma_1(dx^3)^2+
2(\epsilon x^1-2\epsilon x^2\Sigma_1)dx^3dx^4+
$$
$$
4((\epsilon  x^2)^2   \Sigma_1+{\epsilon}^2
x^1x^2-\frac{3}{2}\epsilon x^1 A \Sigma_1)(dx^4)^2\rbrace+3Adx^3dx^4+
12\epsilon x^2  A(dx^4)^2+\sum_{\sigma}e_{\sigma} 
{\Pi'}_{i}(f_i-f_\sigma)(dx^{\sigma})^2,
$$
where
\bqq{2.46}
A=\epsilon x^3+\theta(x^4),
\quad
\Sigma_1=(f_5-\epsilon    x^4)^{-1}+(f_6-\epsilon   x^4)^{-1},
\eeq
$f_1=f_2=f_3=f_4=\epsilon x^4$, $f_5$ is a function of the variable $x^5$,
$f_6$ is a function of the variable $x^6$,  $e_4, e_5, e_6=\pm 1$,
$\epsilon$ equal to 0 or 1.

In this case, nonzero components of the curvature tensor have the following
form: $R_{\nu i \mu}^i$ $(i \ne
\mu)$, $R_{\nu b_1 \mu}^{a_1}$ $(b_1 \ne \mu, a_1 <
b_1)$. Here $i=1,\ldots,6$, $a_1, b_1=1,\ldots,4$, $\nu, \mu$
are the indices of nonzero components of metric $g_{\nu \mu}$ (\ref{2.45}).

Formula (\ref{4.1}) leads to
\bqq{4.23}
R_{114}^1=R_{1\sigma 4}^{\sigma},
\quad
R_{214}^1=R_{2\sigma 4}^{\sigma},
\quad
R_{224}^1=0.
\eeq
For the $h$-space of the $[411]$ type,
$$
R_{114}^1=\rho_4 g_{14},
\quad
R_{214}^1=\rho_4 g_{24},
\quad
R_{1\sigma 4}^{\sigma}=\rho_{\sigma 4} g_{14},
\quad
R_{2\sigma 4}^{\sigma}=\rho_{\sigma 4} g_{24},
$$
$$
R_{224}^1=\gamma_1 g_{24}+\gamma_2
g_{14}+\frac{2\epsilon}{3A^2}(\theta'-\frac{4}{3}{\epsilon}^2 x^2),
$$
where
\bqq{4.23.1}
\gamma_1=-\frac{1}{4} \sum_{\sigma} \frac{(f'_{\sigma})^2}
{(f_{\sigma}-f_4)^3 g_{\sigma \sigma}},
\quad
\gamma_2=-\frac{1}{4} \sum_{\sigma} \frac{(f'_{\sigma})^2}
{(f_{\sigma}-f_4)^4 g_{\sigma \sigma}},
\eeq
$\rho_4$ and $\rho_{\sigma 4}$ are determined by (\ref{4.4}) when $p=4$,
$\sigma=5, 6$.
Therefore, from first two equalities of (\ref{4.23}), we obtain
$$
\rho_4-\rho_{\sigma 4}=0,
\quad
\gamma_1 g_{14}+\frac{2{\epsilon}^2}{3A}=0.
$$
Differentiating the last equality with respect to $x^1$, one obtains
$\epsilon=0$, hence,
$\gamma_1=0$. Then, from the last equality of (\ref{4.23}) it follows that
$\gamma_2=0$.

Let us prove that the necessary conditions of constant curvature of the
$h$-space of the $[411]$ type,
\bqq{4.24}
\rho_4-{\rho}_{\sigma 4}=\epsilon=\gamma_1=\gamma_2=0
\eeq
are sufficient. One has
$$
\pl_1 \rho_4=\pl_2 \rho_4=\pl_3 \rho_4=0,
$$
$$
\pl_4 \rho_4=6 \epsilon \gamma_1,
\quad
\pl_{\sigma} \rho_4=-f'_{\sigma}\frac{\rho_4-\rho_{\sigma
4}}{f_{\sigma}-f_4}.
$$
 From these relations and taking into account (\ref{4.24}),
one obtains $\rho_4=\rho_{\sigma 4}={\rm
const}$. Using the last equality and the conditions
$\epsilon=\gamma_1=\gamma_2=0$, we can prove that
the curvature tensor of our space satisfies identities (\ref{4.1})
with $K=\rho_4={\rm const}$.

\section{The $h$-space of the $[51]$ type}

Metric of the  h-space of the $[51]$ type is defined by the formulas
\bqq{2.53}
g_{ij}dx^idx^j=e(f_6-\epsilon x^5)
\lbrace 8 Adx^1dx^5+2 dx^2dx^4+
\eeq
$$
2(3\epsilon x^3-4A\Sigma_1)dx^2dx^5 +
(dx^3)^2-2\Sigma_1dx^3dx^4+2(2\epsilon x^2-3\epsilon x^3\Sigma_1)dx^3dx^5+
$$
$$
2(\epsilon x^1-2\epsilon x^2 \Sigma_1)dx^4dx^5+
4(3/2 \epsilon x^1\epsilon x^3+(\epsilon  x^2)^2-
2\epsilon x^1 A \Sigma_1-
$$
$$
3\epsilon x^2\epsilon x^3\Sigma_1)(dx^5)^2\rbrace+
e_6(\epsilon x^5-f_6)^5(dx^6)^2,
$$
where
\bqq{2.54}
A=\epsilon
x^4+\theta,
\quad
\Sigma_1=(f_6-\epsilon    x^5)^{-1},
\eeq
$\theta$ is a function of the variable $x^5$,
$f_5=\epsilon x^5$,   $f_6$  is a function of the variable $x^6$,
$e_5, e_6=\pm 1$,
$\epsilon$ equal to 0 or 1.

The nonzero components of the curvature tensor have the form
$R_{\nu i \mu}^i$ $(i \ne
\mu)$, $R_{\nu b_1 \mu}^{a_1}$ $(b_1 \ne \mu, a_1 <
b_1)$, where $i=1,\ldots,6$, $a_1, b_1=1,\ldots,5$, $\nu, \mu$
are the indices of nonzero components of metric $g_{\nu \mu}$ (\ref{2.53}).
Similarly to the previous cases, one can prove that the necessary and 
sufficient
condition of constant curvature of the  $h$-space of the $[51]$ type is
determined by relations
$$
f'_6=\epsilon=0.
$$

\bigskip

\bigskip

Thus, we come to the following

{\bf Theorem.} {\it
The 6-dimensional rigid $h$-spaces $V^6$ of the $[2211]$,
$[321]$, $[33]$, $[411]$ and $[51]$ types have the constant curvature
if and only if

\noindent
for the $h$-space of the $[2211]$ type
\bqq{4.47}
\rho_p-\rho_{\sigma p}=\rho_p-\rho_{pq}=\epsilon=\tilde \epsilon=0
\quad
(p\ne q, p, q=2, 4, \sigma=5, 6),
\eeq

\noindent
for the $h$-space of the $[321]$ type
\bqq{4.47.1}
{f'}_6=\epsilon=\tilde\epsilon=0,
\eeq

\noindent
for the $h$-space of the $[33]$ type
\bqq{4.47.2}
\epsilon=\tilde\epsilon=0,
\eeq

\noindent
for the $h$-space of the $[411]$ type
\bqq{4.46}
\rho_4-\rho_{\sigma 4}=\epsilon=\gamma_1=\gamma_2=0
\quad
(\sigma=5, 6),
\eeq

\noindent
for the $h$-space of the $[51]$ type
\bqq{4.47.3}
{f'}_6=\epsilon=0,
\eeq
where $\rho_p$, $\rho_{\sigma p}$, $\rho_{pq}$, $\gamma_1$ and $\gamma_2$
are determined by formulas {\rm (\ref{4.4})}, {\rm(\ref{4.23.1})}.
}

\bigskip

I am grateful to A.V.Aminova for the constant encouragement and
discussions.

\end{document}